\documentclass[11pt]{amsart}
\usepackage{amssymb}
\usepackage{bm}
\usepackage{graphicx}
\usepackage[centertags]{amsmath}
\usepackage{amsfonts}
\usepackage{amsthm}

\theoremstyle{plain}% default
\theoremstyle{definition}
\newtheorem{theorem}{Theorem}[section]
\newtheorem{lemma}[theorem]{Lemma}
\newtheorem{proposition}[theorem]{Proposition}
\newtheorem{corollary}[theorem]{Corollary}
\newtheorem{definition}[theorem]{Definition}
\newtheorem{remark}[theorem]{Remark}

\newtheorem{notation}[theorem]{Notation}
\newtheorem*{claim}{Claim}

\newcommand{\e}{\varepsilon}
\newcommand{\bb}{\overline{b}}

\DeclareMathOperator{\A}{\mathcal{A}}
\DeclareMathOperator{\F}{\mathcal{F}}
\DeclareMathOperator{\T}{\mathcal{T}}
\DeclareMathOperator{\N}{\mathbb{N}}

\DeclareMathOperator{\R}{\mathbb{R}}

\DeclareMathOperator{\maxsupp}{maxsupp}
\DeclareMathOperator{\minsupp}{minsupp}
\DeclareMathOperator{\supp}{supp}

\DeclareMathOperator{\ran}{ran}
\DeclareMathOperator{\rank}{rank}

\newcommand{\bp}{\begin{proof}}
\newcommand{\ep}{\end{proof}}
\newcommand{\bd}{\begin{definition}}
\newcommand{\ed}{\end{definition}}
\newcommand{\bt}{\begin{theorem}}
\newcommand{\et}{\end{theorem}}
\newcommand{\bpr}{\begin{proposition}}
\newcommand{\epr}{\end{proposition}}
\newcommand{\bc}{\begin{corollary}}
\newcommand{\ec}{\end{corollary}}
\newcommand{\bl}{\begin{lemma}}
\newcommand{\el}{\end{lemma}}
\newcommand{\br}{\begin{remark}}
\newcommand{\er}{\end{remark}}
\newcommand{\be}{\begin{enumerate}}
\newcommand{\ee}{\end{enumerate}}

\begin{document}

\title{More $\ell_r$ saturated $\mathcal{L}^{\infty}$ spaces}
\author{I. Gasparis}
\address{Department of Mathematics, Aristotle University of Thessaloniki, Thessaloniki 54124, Greece} 
\email{ioagaspa@math.auth.gr}
\author{M.K. Papadiamantis}
\address{National Technical University of Athens, Faculty of Applied Sciences,
Department of Mathematics, Zografou Campus, 157 80, Athens,
Greece}
\email{mpapadiamantis@yahoo.gr}
\author{D.Z. Zisimopoulou}
\address{National Technical University of Athens, Faculty of Applied Sciences,
Department of Mathematics, Zografou Campus, 157 80, Athens,
Greece}
\email{dzisimopoulou@hotmail.com}


%------------------------Abstract-------------------------------%

\begin{abstract} Given $r \in (1, \infty)$, we construct a new $\mathcal{L}^{\infty}$ separable
Banach space which is
$\ell_r$ saturated .
\end{abstract}

\maketitle

%-------------------------Introduction---------------------------%

\section{Introduction}
The Bourgain-Delbaen spaces \cite{BD} are examples of separable $\mathcal{L}^{\infty}$ spaces
containing no isomorphic copy of $c_0$. They have played a key role
in the solution of the scalar-plus-compact problem by Argyros and Haydon \cite{AH},
where a Hereditarily
Indecomposable $\mathcal{L}^{\infty}$ space is presented with the property that every
operator on the space is a compact perturbation of a scalar multiple of the identity.

There has recently been an interest in the study $\mathcal{L}^{\infty}$ spaces
of the Bourgain-Delbaen type. Freeman, Odell and Schlumprecht \cite{fos} showed that
every Banach space with separable dual is isomorphic to a subspace of a
$\mathcal{L}^{\infty}$ space having a separable dual.
The aim of this paper is to present a method of constructing, for every $1<r<\infty$, a
new $\mathcal{L}^{\infty}$ space which is $\ell_r$ saturated.
Our approach shares common features with the
Argyros-Haydon work. More
precisely we combine, as in \cite{AH}, the Bourgain-Delbaen method
\cite{BD} yielding exotic $\mathcal{L}^{\infty}$ spaces, with the
Tsirelson type norms that are equivalent to some $\ell_r$ norm
(see \cite{AD}, \cite{AT}, \cite{B}). Recall that in \cite{H},
the original Bourgain-Delbaen spaces $\mathfrak{X}_{a,b}$ with
$a<1$, $b<\frac{1}{2}$ and $a+b>1$ where shown to be $\ell_p$ saturated for
$p$ determined by the formulas $\frac{1}{p}+\frac{1}{q}=1$ and
$a^q+b^q=1$.

This paper is organized as follows. In the second section, for a
given $r\in(1,\infty)$, we construct a Banach space
$\mathfrak{X}_r$. To do this, we first choose $n\in\N$,
$n>1$, and a finite sequence $\bb=(b_1,b_2,\ldots,b_n)$ of positive
real numbers with $b_1<1$, $b_2,b_3,\ldots,b_n<\frac{1}{2}$
such that $\sum_{i=1}^nb_i^{r'}=1$ and
$\frac{1}{r}+\frac{1}{r'}=1$. The definition of
$\mathfrak{X}_r$ combines the Bourgain-Delbaen method with the
Tsirelson type space $\T(\A_n,\bb)$ which will be later proved to be
isomorphic to $\ell_r$. In particular, if
$b_1=b_2=\ldots=b_n=\theta$, $\T(\A_n,\bb)$ coincides
with $\T(\A_n,\theta)$ and the latter is known to be isomorphic to $\ell_p$ for some
$p\in(1,\infty)$ (see \cite{AT}). It is worth
noticing that for $n=2$ the spaces $\mathfrak{X}_r$ essentially coincide with
the original Bourgain-Delbaen spaces $\mathfrak{X}_{a,b}$. Thus,
our construction of
$\mathcal{L}^{\infty}$ spaces which are $\ell_r$ saturated spaces, can
be considered as a generalization
of the Bourgain-Delbaen method. We must point out here that when $n=2$, our
proof of the fact that $\mathfrak{X}_r$ is $\ell_r$ saturated, differs from
Haydon's (see \cite{H}) corresponding one for $\mathfrak{X}_{a,b}$.
To be more specific, $\mathfrak{X}_r$ has a natural FDD $(M_k)$. Given a normalized
skipped block basis $(u_k)$ of $(M_k)$ with the supports of the $u_k$'s lying far enough apart,
then it is not hard to check that $(u_k)$ dominates $(e_k)$, the natural basis of $\T(\A_n,\bb)$. The
same holds for every normalized block basis of $(u_k)$. To obtain a normalized block basis of $(u_k)$
equivalent to $(e_k)$, we
select a sequence $I_1 < I_2 < \dots $ of successive finite subsets of $\mathbb{N}$ such that
$\lim_k \|\sum_{i \in I_k} u_i \| = \infty$. Such a choice is possible by the domination of
$(e_k)$ by $(u_k)$. We set $v_k = \|\sum_{i \in I_k} u_i \|^{-1} \sum_{i \in I_k} u_i $
and show that some subsequence of $(v_k)$ is dominated by $(e_k)$.
To accomplish this we adapt the method of the analysis of the members of a finite
block basis of $(e_k)$ with respect to a
functional in the natural norming
set of $\T(\A_n,\bb)$ (see \cite{DB}), to the context of the present construction.
We believe that this approach yields a more transparent proof than Haydon's, at least for the upper
$\ell_r$ estimate.

The rest of the paper is devoted to the proof of the main
property, namely that $\mathfrak{X}_r$ is $\ell_r$ saturated. In
Section 3, we define the tree analysis of the functionals
$\{e_{\gamma}^*:\gamma\in\Gamma\}$ which is a 1-norming subset of
the unit ball of $\mathfrak{X}_r^*$. The tree analysis is similar
to the corresponding one used in the Tsirelson
and mixed Tsirelson spaces \cite{AT}. In the following two
sections we establish the lower and upper norm estimates for
certain block sequences in the space $\mathfrak{X}_r$.

In the final section we show that every block basis of $(M_k) $ admits
a further normalized block basis
$(x_k)$ such that every normalized block basis of $(x_k)$ is
equivalent to the natural basis of the space $\T(\A_n,\bb)$. Zippin's
theorem \cite{Z} yields the desired result.

%-------------------------Definition Of The Space---------------------------%

\section{Preliminaries}

In this section we define the space $\mathfrak{X}_r$ combining the
Bourgain-Delbaen construction \cite{BD} and the Tsirelson type
constructions \cite{AD}, \cite{AT}.

Before proceeding, we recall some notation and terminology from \cite{AH}. Let
$n\in\N$ and $0<b_1,b_2,...,b_n<1$ with $\sum_{i=1}^nb_i>1$ and there exists
$r'\in(1,\infty)$ such that $\sum_{i=1}^n{b_i}^{r'}=1$. We may also assume without loss of generality that $b_1>b_2>\ldots>b_n$.
We define $W[(\A_n,\bb)]$
to be the smallest subset $W$ of $c_{00}(\N)$ with the following properties:

\begin{enumerate}
\item $\pm e^*_k\in W$ for all $k\in \mathbb N$, \item whenever
$f_i\in W$ and $\max \supp f_i<\min \supp f_{i+1}$ for all $i$, we
have $\sum_{i\le a} b_if_i \in W$, provided that $a\leq n$,
\end{enumerate}

We say that an element $f$ of $W[(\A_n,\bb)]$ is of Type $0$ if $f=\pm e_k^*$ for some $k$ and of Type I otherwise; an element of Type I is said to have weight $b_a$ for some $a\leq n$ if $f=\sum_{i=1}^af_i$ for a suitable sequence $(f_i)$ of successive elements of $W[\A_n,\bb]$.\\
The {\em Tsirelson space} $\T(\A_n,\bb)$ is defined to be the
completion of $c_{00}$ with respect to the norm
\[\|x\|= \sup \{
\langle f,x\rangle: f\in W[\A_n,\bb]\}.\] We may also characterize
the norm of this space implicitly as being the smallest function
$x\mapsto \|x\|$ satisfying
\[\|x\| =
\max\bigg\{\|x\|_\infty,\sup\sum_{i=1}^n b_i\|E_ix\|\bigg\},\]
where the supremum is taken over all sequences of finite subsets
$E_1<E_2<\cdots<E_n$.

We shall now present the fundamental aspects related
to the Bourgain-Delbaen construction.\\
For the interested readers we mention that the following method
can be characterized as the "dual" construction of the
construction presented in \cite{AH}. This characterization is
based on the fact that in \cite{AH} a particular kind of basis is
given to $\ell_1(\Gamma)$ and the Bourgain-Delbaen type space $X$
is seen as the predual of its dual, which is $\ell_1(\Gamma)$.

Let $(\Gamma_q)_{q\in{\N}}$ be a strictly increasing sequence of
finite sets and denote their union by $\Gamma$;
$\Gamma=\cup_{q\in{\N}} \Gamma_q$.\\
We set $\Delta_0=\Gamma_0$ and
$\Delta_q=\Gamma_q\backslash\Gamma_{q-1}$ for $q=1,2,\ldots$\\
Assume furthermore that to each $\gamma\in\Delta_q$, $q\geq1$, we
have assigned a linear functional
$c_{\gamma}^*:\ell^{\infty}(\Gamma_{q-1})\rightarrow\R$. Next, for
$n<m$ in $\N$, we define by induction, a linear operator
$i_{n,m}:\ell^{\infty}(\Gamma_{n})\rightarrow\ell^{\infty}(\Gamma_{m})$
as follows:\\
For $m=n+1$, we define
$i_{n,n+1}:\ell^{\infty}(\Gamma_{n})\rightarrow\ell^{\infty}(\Gamma_{n+1})$
by the rule
\[(i_{n,n+1}(x))(\gamma)=
\begin{cases}
x(\gamma),\ \text{if }\gamma\in\Gamma_n\\
c_{\gamma}^*(x),\ \text{if }\gamma\in\Delta_{n+1}
\end{cases}\]
for every $x\in\ell^{\infty}(\Gamma_n)$.\\
Then assuming that $i_{n,m}$ has been defined, we set
$i_{n,m+1}=i_{m,m+1}\circ i_{n,m}$. A direct consequence of the
above definition is that for $n<l<m$ it holds that
$i_{n,m}=i_{l,m}\circ i_{n,l}$. Finally we denote by
$i_n:\ell^{\infty}(\Gamma_n)\rightarrow\R^{\Gamma}$ the direct
limit $i_n=\lim_{m\rightarrow\infty}i_{n,m}$.

We assume that there exists a $C>0$ such that for every $n<m$ we
have $\|i_{n,m}\|\le C$. This implies that $\|i_n\|\le C$
and therefore
$i_n:\ell^{\infty}(\Gamma_n)\to\ell^{\infty}(\Gamma)$ is a bounded
linear map. In particular, setting
$X_n=i_n[\ell^{\infty}(\Gamma_n)]$, we have that
$X_n\stackrel{C}\thickapprox\ell^{\infty}(\Gamma_n)$ and
furthermore $(X_n)_{n\in\N}$ is an increasing sequence of
subspaces of $\ell^{\infty}(\Gamma)$. We also set
$\mathfrak{X}_{BD}=\overline{\bigcup\limits_{n\in\N}X_n}\hookrightarrow\ell^{\infty}(\Gamma)$
equipped with the supremum norm. Evidently,
$\mathfrak{X}_{BD}$ is an $\mathcal{L}^{\infty}$ space.

Let us denote by
  $r_n:\ell^{\infty}(\Gamma)\to\ell^{\infty}(\Gamma_n)$ the natural
  restriction map, i.e. $r_n(x)=x|_{\Gamma_n}$. We will also abuse
  notation and denote by
  $r_n:\ell^{\infty}(\Gamma_m)\to\ell^{\infty}(\Gamma_n)$ the
  restriction function from $\ell^{\infty}(\Gamma_m)$ to $\ell^{\infty}(\Gamma_n)$
  for $n<m$.

  \begin{notation}
  \item[(i)] We denote by $e_\gamma^*$ the restriction of the unit
  vector $e_\gamma\in\ell^1(\Gamma)$ on the space $\mathfrak{X}_{BD}$.
  \item[(ii)] We also extend the functional
  $c_\gamma^*:\ell^{\infty}(\Gamma_n)\to\R$ to a functional
  $c_\gamma^*:\mathfrak{X}_{BD}\to\R$ by the rule
  $c_\gamma^*(x)=(c_\gamma^*\circ r_{q-1})(x)$ when $\gamma\in
  \Delta_q$.
  \end{notation}

As it is well known from \cite{AH} and \cite{BD}, instead of the
Schauder basis of $\mathfrak{X}_{BD}$, it is more convenient to
work with a FDD naturally defined as follows:

For each $q\in\N$ we set $M_q=i_q[\ell^{\infty}(\Delta_q)]$.\\
We
briefly establish this fact in the following proposition and then
continue with the details of the construction of $\mathfrak{X}_r$.

\bpr The sequence $(M_q)_{q\in\N}$ is a FDD for
$\mathfrak{X}_{BD}$. \epr

\bp For $q\geq 0$ we define the maps $P_{\{q\}}:\mathfrak{X}_{BD}
\to M_q$ with
 \[P_{\{q\}}(x)=i_q(r_q(x))-i_{q-1}(r_{q-1}(x))\]

 It is easy to check that each $P_{\{q\}}$ is a projection onto $M_q$ and that for $q_1\neq q_2$ and $x\in M_{q_2}$ we have
 $P_{\{q_1\}}(x)=0$. Also we have that $\|P_{q}\|\leq 2 C$. We point out that in a similar manner one can define projections on intervals of the form $I=(p,q]$ so that
 $P_I(x)=\sum_{i=p+1}^q P_{\{i\}}(x)$ for which we can readily verify the formula
 \[P_I(x)=i_q(r_q(x))-i_p(r_p(x))\]
 Note that $\|P_I\|\leq 2C$. This shows that indeed $(M_q)_q$ is a FDD generating $\mathfrak{X}_{BD}$.
 \ep

For $x\in\mathfrak{X}_{BD}$ we denote by $\supp x$ the set $\supp
x=\{q:P_{\{q\}}(x)\neq0\}$ and by $\ran x$ the minimal interval of
$\N$ containing $\supp x$.

\begin{definition}
A block sequence $(x_i)_{i=1}^{\infty}$ in $\mathfrak{X}_{BD}$ is
called $\emph{skipped}$ (with respect to $(M_q)_{q\in\N}$), if
there is a subsequence $(q_i)_{i=1}^{\infty}$ of $\N$ so that for
all $i\in\N$, $\maxsupp{x_i}<q_i<\minsupp{x_{i+1}}$.
\end{definition}
In the sequel, when we refer to a skipped block sequence, we
consider it to be with respect to the FDD $(M_q)_{q\in\N}$.

Let $q\geq 0$. For all $\gamma\in\Delta_q$ we set
$d_\gamma^*=e_{\gamma}\circ P_{\{q\}}.$ Then the family
$(d^*_{\gamma})_{\gamma\in \Gamma}$ consists of the biorthogonal
functionals of the FDD $(M_q)_{q\geq 0}$. Notice that for
$\gamma\in\Delta_q$,

\begin{eqnarray*}
d_\gamma^*(x)& = & P_q(x)(\gamma)=i_q(r_q(x))(\gamma)-i_{q-1}(r_{q-1}(x))(\gamma)=\\
& =& r_q(x)(\gamma)-c_\gamma^*(r_{q-1}(x))=x(\gamma)-c_\gamma^*(x)=\\
& = & e^*_{\gamma}(x)-c^*_{\gamma}(x).
\end{eqnarray*}

The sequences $(\Delta_q)_{q\in\N}$ and
$(c_{\gamma}^*)_{\gamma\in\Gamma}$ are determined as in \cite{AH}, section 4 and Theorem 3.5.\\
We give some useful notation. For fixed $n\in\N$ and
$\bb=(b_1,b_2,\ldots,b_n)$ with $0<b_1,b_2,\ldots,b_n<1$, for each
$\gamma\in \Delta_q$ we assign
\begin{enumerate}
     \item[(a)] $\rank \gamma=q$
     \item[(b)] age of $\gamma$ denoted by $a(\gamma)=a$ such that
     $2\leq a\leq n$
     \item[(c)] weight of $\gamma$ denoted by $w(\gamma)=b_a$
\end{enumerate}
In order to proceed to the construction, we first need to fix a
positive integer $n$ and a descending sequence of positive real numbers $b_1,\ldots,b_n$
such that $b_1<1$, $b_i<\frac{1}{2}$, for every $i=2,\ldots,n$ and
$\sum_{i=1}^nb_i>1$. Let $r\in(1,\infty)$ be such that
$\sum_{i=1}^nb_i^{r'}=1$ and $\frac{1}{r}+\frac{1}{r'}=1$. Now we
shall define the space $\mathfrak{X}_r$ by using the Bourgain-Delbaen
construction that was presented in the preceding paragraphs.

We set $\Delta_0=\emptyset$, $\Delta_1=\{0\}$ and recursively
define for each $q>1$ the set $\Delta_q$.\\
Assume that $\Delta_p$ have been defined for all $p\le q$. We set
\begin{eqnarray*}
     \Delta_{q+1} & = & \big\{ (q+1,a,p,\eta,\e e_{\xi}^*):\ 2\leq a\leq n, p\leq q,\ \e=\pm1,\ e_{\xi}^*\in S_{\ell^1(\Gamma_q)},\ \xi\in\Gamma_q\setminus\Gamma_p,\\
      & & \ \ \eta\in\Gamma_p,\ b_
      {a-1}=w(\eta)\big\}
 \end{eqnarray*}
For $\gamma\in\Delta_{q+1}$ it is clear that the first coordinate
is the $\rank$ of $\gamma$, while the second is the age
$a(\gamma)$ of $\gamma$. The functionals
$(c_{\gamma}^*)_{\gamma\in\Delta_{q+1}}$ are defined in a way that
depends on $\gamma\in \Delta_{q+1}$. Namely, let $x\in
\ell^{\infty}(\Gamma_q)$.
        \begin{enumerate}
        \item[(i)] For $\gamma=(q+1,2,p,\eta,\e e_{\xi}^*)$ we set
        \[
        c_\gamma^*(x)=b_1x(\eta)+b_2\e e_{\xi}^*\big( x-i_{p,q}(r_p(x))\big).\]
         \item[(ii)] For $\gamma=
             (q+1,a,p,\eta,\e e_{\xi}^*)$ with $a>2$ we set
          \[
        c_\gamma^*(x)=x(\eta)+b_a\e e_{\xi}^*\big(x-i_{p,q}(r_p(x))\big).\]
        \end{enumerate}

We may now define sequences $(i_q)$, $(\Gamma_q)$, $(X_q)$ in
a similar manner as before and set
$\mathfrak{X}_r=\overline{\bigcup\limits_{q\in\N}X_q}$. Assuming
that $(i_q)$ is uniformly bounded by a constant C, we conclude
that the space $\mathfrak{X}_r$ is a subspace of
$\ell_{\infty}(\Gamma)$. The constant C is determined as in
\cite{AH} Theorem 3.4, by taking $C=\frac{1}{1-2b_2}$. Thus, for
every $m\in\N$, $\|i_m\|\leq C$. This implies that $\|P_I\|\leq
2C$ for every $I$ interval.

\br In the case of $n=2$, i.e. $\bb=(b_1,b_2)$, the space
$\mathfrak{X_r}$ essentially coincides with the Bourgain-Delbaen space
$\mathfrak{X}_{b_1,b_2}$, since every $\gamma\in\Gamma$ is of age
2. \er

\br As it is shown in Proposition \ref{13}, the choice of r, based
on the fixed $n$ and $\bb$, yields that
$\T(\mathcal{A}_n,\bb)\cong\ell_r$. Moreover, the ingredients of
the "Tsirelson type spaces" theory that are used throughout this
paper are essentially the same with the corresponding ones in \cite{AH}. The
basic difference in our approach is that we use only one family
$\T(\A_n,\bb)$ for some appropriate $n$ and $\bb$. \er

%-------------------------Tree Analysis---------------------------%

\section {The Tree Analysis of $e_{\gamma}^*$ for $\gamma\in\Gamma$}
We begin by recalling the analysis of $e_{\gamma}^*$ in \cite{AH}
section 4. The only difference is that in our case all the
functionals $e_{\gamma}^*$ have weight depending on their age which
is greater or equal to 2.
\subsection{The evaluation Analysis of $e_{\gamma}^*$ for $\gamma\in\Gamma$} First we point out that for $q\in\N$ every
$\gamma\in\Delta_{q+1}$ admits a unique analysis as follows:\\
Let $a(\gamma)=a\leq{n}$. Then using backwards induction we
determine a sequence of sets
$\{p_i,q_i,\varepsilon_ie_{\xi_i}^*\}_{i=1}^a\cup\{\eta_i\}_{i=2}^a$
with the following properties.
    \begin{enumerate}
    \item[(i)] $p_1<q_1<\cdots<p_a<q_a=q$.
    \item[(ii)] $\varepsilon_i=\pm1,\ \rank \xi_i\in (p_i,q_i]$ for $1\le i\le
    a$ and $\rank \eta_i=q_i+1$ for $2\le i\le a$.
    \item[(iii)] $\eta_a=\gamma$, $\eta_i=(\rank\eta_i,i,p_i,\eta_{i-1},\varepsilon_ie_{\xi_i}^*)$ for every
    $i>2$\\ $\eta_2=(\rank\eta_2,2,p_2,\e_1\xi_1,\varepsilon_2e_{\xi_2}^*)$ and $(p_1,q_1]=(1,\rank\xi_1]$.
    \end{enumerate}

\begin{definition}
Let $q\in\N$ and $\gamma\in\Gamma_q$. Then the sequence
$\{p_i,q_i,\varepsilon_ie_{\xi_i}^*\}_{i=1}^a\cup\{\eta_i\}_{i=2}^a$
satisfying all the above properties will be called the analysis of
$\gamma$.
\end{definition}

Moreover, following similar arguments as in \cite{AH} Proposition 4.6 it holds that,
\[e_\gamma^*=\sum_{i=2}^a d_{\eta_i}^*+\sum_{i=1}^ab_i\varepsilon_ie_{\xi_i}^*\circ
P_{(p_i,q_i]}=\sum_{i=2}^a e_{\eta_i}^*\circ P_{\{q_i+1\}}
+\sum_{i=1}^ab_i\varepsilon_ie_{\xi_i}^*\circ P_{(p_i,q_i]}.\]
We set $g_{\gamma}=\sum\limits_{i=2}^a d_{\eta_i}^*$ and $f_{\gamma}=\sum\limits_{i=1}^ab_i\varepsilon_ie_{\xi_i}^*\circ P_{(p_i,q_i]}$.\\

\subsection{The r-Analysis of the functional $e_{\gamma}^*$}
Let $r\in{\N}$ and $\gamma\in\Delta_{q+1}$.Let
$a(\gamma)=a\leq{n}$ and
$\{p_i,q_i,\varepsilon_ie_{\xi_i}^*\}_{i=1}^a\bigcup\{\eta_i\}_{i=2}^a$
the evaluation analysis of $\gamma$. We define the r-analysis of
$e_{\gamma}^*$ as follows:
\begin{enumerate}
\item[(a)] If $r\leq{p_1}$, then the r-analysis of $e_{\gamma}^*$
coincides with the evaluation analysis of $e_{\gamma}^*$.
\item[(b)] If $r\geq{q_a}$, then we assign no r-analysis to
$e_{\gamma}^*$ and we say that $e_{\gamma}^*$ is r-indecomposable.
\item[(c)] If $p_1<r<q_a$, we define $i_r=\min\{i:r<q_i\}$. Note that this
is well-defined. The r-analysis of $e_{\gamma}^*$ is the following
triplet
\[\{(p_i,q_i]\}_{i\geq i_r},\{\e_i\xi_i\}_{i\geq i_r},\{\eta_i\}_{i\geq\max\{2,i_r\}}.\]
where $p_{i_r}$ is either the same or $r$ in the case that
$r>p_{i_r}$.
\end{enumerate}

Next we introduce the tree analysis of $e_{\gamma}^*$ which is
similar to the tree analysis of a functional in a Mixed Tsirelson
space (see \cite{AT} Chapter II.1). Notice that the evaluation
analysis and the r-analysis of $e_{\gamma}^*$ form the first
level of the
tree analysis that we are about to present. \\

We start with some notation. We denote by $(\T,"\preceq")$ a
finite partially ordered set which is a tree. Its elements are finite sequences of natural numbers ordered by
the initial segment partial order. For every
$t\in\T$,we denote by $S_t$ the immediate successors of $t$

Assume now that
$(p_t,q_t]_{t\in\T}$ is a tree of intervals of $\N$ such that
$t\preceq s$ iff $(p_t,q_t]\supset(p_s,q_s]$ and $t,s$ are
incomparable iff $(p_t,q_t]\cap(p_s,q_s]=\emptyset$. For such a
family $(p_t,q_t]_{t\in\T}$ and $t,s$ incomparable we shall denote
by $t<s$ iff $(p_t,q_t]<(p_s,q_s]$ (i.e. $q_t<p_s$).

\subsection{The Tree Analysis of the functional $e_{\gamma}^*$}
Let $\gamma\in\Delta_{q+1}$ with $a(\gamma)=a\leq{n}$. A family of
the form
$\mathcal{F}_{\gamma}=\{\xi_t, (p_t,q_t]\}_{t\in\T}$ is called the tree analysis of $e_{\gamma}^*$ if the following are satisfied:\\
\begin{enumerate}
\item $\T$ is a finite tree with a unique root denoted as
$\emptyset$. \item We set
$\xi_{\emptyset}=\gamma$,$(p_{\emptyset},q_{\emptyset}]=(1,q]$ and
let
$\{p_i,q_i,\varepsilon_ie_{\xi_i}^*\}_{i=1}^a\bigcup\{\eta_i\}_{i=2}^a$
the evaluation analysis of $\xi_{\emptyset}$. Set
$S_{\emptyset}=\{(1),(2),\ldots,(a)\}$ and for every $s=(i)\in
S_{\emptyset}$, $\{\xi_s,(p_s,q_s]\}=\{\xi_i,(p_i,q_i]\}$.

\item Assume that for a $t\in\T$ $\{\xi_t,(p_t,q_t]\}$ has been
defined. There are two cases:\\
\begin{enumerate}
\item If $e_{\xi_t}^*$ is $p_t$-decomposable, let
\[\{(p_i,q_i]\}_{i\geq i_{p_t}},\{\e_i\xi_i\}_{i\geq
i_{p_t}},\{\eta_i\}_{i\geq\max\{2,i_{p_t}\}}\] the $p_t$ analysis
of $e_{\xi_t}^*$. We set $S_t=\{(t^\smallfrown i): i\geq
i_{p_t}\}$ and
\[S_t^{p_t}=\begin{cases} S_t,\ \text{if }
\eta_{i_{p_t}}\ \text{exists}\\
S_t\diagdown\{(t^\smallfrown i_{p_t})\},\ \text{otherwise}
\end{cases}\]
Then, for every $s=(t^\smallfrown i)\in S_t$, we set
$\{\xi_s,(p_s,q_s]\}=\{\xi_i,(p_i,q_i]\}$ where
$\{\e_i\xi_i,(p_i,q_i]\}$ is a member of the $p_t$ analysis of
$e_{\xi_t}^*$.
 \item $e_{\xi_t}^*$ is $p_t$-indecomposable, then
$\xi_t$ consists a maximal node of $\F_{\gamma}$.
\end{enumerate}
\end{enumerate}

\begin{notation} For later use we need the
following:\\
For every $t\in\T$ $e_{\xi_t}^*=f_t+g_t$, where $f_t=\sum_{s\in
S_t}b_s\e_se_{\xi_s}^*\circ P_{(p_s,q_s]}$ and $g_t=\sum_{s\in
S_t^{p_t}}d_{\eta_s}^*$ and for
$s=(t^\smallfrown i)\in S_t^{p_t}$, \\
$\eta_{(t^\smallfrown i)}=(\rank\eta_{(t^\smallfrown
i)},i,p_{(t^\smallfrown i)},\eta_{(t^\smallfrown
i-1)},\e_{(t^\smallfrown i)}e_{\xi_{(t^\smallfrown i)}}^*)$.\\ In
the rest of the paper, we set $f_t=f_{\xi_t}$ and $g_t=g_t$.
\end{notation}

\bl \label{4} Let $x\in\mathfrak{X}_r$ and $\gamma\in\Gamma$.
Then,
\[e_{\gamma}^*(x)=\prod_{\emptyset\preceq s\preceq
t_x}(\e_sb_s)(f_{t_x}+g_{t_x})(x),\] where $t_x=\max\{t:\ran
x\subseteq(p_t,q_t]\}$. \el

\bp Let $\F_{\gamma}=\{\xi_t, (p_t,q_t]\}_{t\in\T}$ a tree
analysis of $\gamma$.\\If $\{t:\ran
x\subseteq(p_t,q_t]\}=\emptyset$, then
$e_{\gamma}^*(x)=f_{\emptyset}(x)+g_{\emptyset}(x)$ and the
equality
holds.\\
If $\{t:\ran x\subseteq(p_t,q_t]\}\neq\emptyset$, we can find
$\{t_1\prec t_2\prec\ldots\prec t_m\}\in\T$ such that $t_1\in
S_\emptyset$ and
$t_m=t_x$.\\
For every $t\in\T$ with $t\prec t_x$, $g_t(x)=0$. Indeed, for
every $s\in S_t^{p_t}$, $d_{\eta_s}^*(x)=e_{\eta_s}^*\circ
P_{\{q_s+1\}}(x)=0$ because
$\ran x\subseteq(p_{t_x},q_{t_x}]\subseteq(p_s,q_s]$.\\
So, we have that
\begin{eqnarray*}
e_{\gamma}^*(x) &=& f_{\emptyset}(x)=\sum_{s\in
S_{\emptyset}}b_s\e_se_{\xi_s}^*\circ P_{(p_s,q_s]}(x)=b_{t_1}
\e_{t_1}e_{\xi_{t_1}}^*(x)\\ &=&b_{t_1}
\e_{t_1}f_{t_1}(x)=b_{t_1}\e_{t_1}b_{t_2}\e_{t_2}e_{\xi_{t_2}}^*\circ
P_{(p_{t_2},q_{t_2}]}(x)= b_{t_1}b_{t_2}\e_{t_1}\e_{t_2}e_{\xi_{t_2}}^*(x)\\
&=&
b_{t_1}b_{t_2}\e_{t_1}\e_{t_2}f_{t_2}(x)=\ldots=\prod_{\emptyset\preceq
s\preceq t_x}(\e_sb_s)(f_{t_x}+g_{t_x})(x)
\end{eqnarray*}
setting $\e_{\emptyset}=b_{\emptyset}=1$. \ep

\bc \label{5} If $(f_{t_x},(p_{t_x},q_{t_x}])$ is a maximal node,
then $e_{\gamma}^*(x)=0$. \ec

\bp Let $(f_{t_x},(p_{t_x},q_{t_x}])$ be a maximal node. Then
$f_{t_x}(x)=0$ and $g_{t_x}(x)=0$ and from Lemma \ref{4} we deduce
that $e_{\gamma}^*(x)=0$. \ep

%-------------------------Lower Inequality---------------------------%

\section{The lower estimate}

\bd An $\phi\in W(\A_n,\bb)$ is said to be a proper functional if
it admits a tree analysis $(\phi_t)_{t\in\T}$ such that for every
non-maximal node $t\in\T$ the set $\{\phi_s:s\in S_t\}$ has at
least two non-zero elements. \ed

We denote by $W_{pr}(\A_n,\bb)$ to be the subset of $W(\A_n,\bb)$
consisting of all proper functionals. For every $t\in\T$ it holds
that $\phi_t=\sum_{s\in S_t} b_s\phi_s$ with $\{b_s\}_{s\in
S_t}\subseteq\{b_1,b_2,\ldots,b_n\}$ and $b_{\emptyset}=1$.

\bl \label{1} The set $W_{pr}(\A_n,\bb)$ 1-norms the space
$\T(\A_n,\bb)$. \el

\bp We shall show that for every $\phi\in W(\A_n,\bb)$ there
exists $g\in W_{pr}(\A_n,\bb)$ such that $|\phi(m)|\leq g(m)\
\forall m\in{\N}$. Since the basis is 1-unconditional the previous statement
yields the result.

To this end, let $\phi\in W(\A_n,\bb)$. Then using a tree analysis
$\{\phi_t\}_{t\in\T}$ of $\phi$ we easily see that for every
$m\in\supp f$, there exists a maximal node $t_m\in\T$ with
$\phi_{t_m}=\e_me_m^*$ and $\phi(m)=\e_m\prod\limits_{t<t_m}b_t$.

For every $m\in\supp\phi$ we set $K_m=\{t\in \T:t<t_m\ \text{and}\
\# S_t>1\}$. Then it is easy to see that the functional
$g=\sum\limits_{m\in\supp\phi}(\prod\limits_{t\in K_m}b_t)e_m^*$
is a functional belonging to $W_{pr}(\A_n,\bb)$. Moreover, since
$b_t<1$ for every $t\in\T$ we get that $|\phi(m)|\leq g(m)\
\forall m\in{\N}$. \ep

\bl \label{2} Let $\phi\in W_{pr}(\A_n,\bb)$ and $l\in\N$. If
$\maxsupp\phi=l$, then $h(\T_{\phi})\leq l$. \el

\bp Let $\theta_n$ be the amount of nodes at the $n$ level of
$\T_{\phi}$. Since $\phi$ is proper, it holds that
$\theta_{n+1}>\theta_n$ for every $n\in\N$. Assume to the contrary
that $h(\T_{\phi})>l$, i.e. $h(\T_{\phi})=l+k$ for some $k\in\N$.
Then,
\[\theta_1=1,\ \theta_2\geq2,\ \ldots\ ,\ \theta_{l+k}\geq l+k\]
Since, the $l+k$ level of $\T_{\phi}$ consists of functionals of
the form $e_i^*$, we deduce that $\maxsupp\phi\geq l+k>l$, which
leads to a contradiction. \ep

\bpr \label{3} Let $(x_k)_{k\in\N}$ be a normalized skipped block
sequence in $\mathfrak{X}_r$ and $(q_k)_{k\in\N}$ a strictly
increasing sequence of integers such that $\supp x_k\subset
(q_k+k,q_{k+1})$. Then, for every sequence of
positive scalars $(a_k)_{k\in\N}$ and for every $l\in\N$, it holds that
\begin{equation}\label{eq1}
\|\sum_{k=1}^la_ke_k\|_{\T(\A_n,\bb)}\leq C\|\sum_{k=1}^la_kx_k\|_{\infty}
\end{equation}
where $(e_k)_{k\in\N}\subseteq\T(\A_n,\bb)$ and C is an upper bound for the norms
of the operators $i_m$ in $\mathfrak{X}_r$.
\epr

\bp Let $\phi\in W(\A_n,\bb)$. From Lemma \ref{1} we may assume
that $\phi$ is proper. We will use induction on the height of the
tree $\T_{\phi}$.

If $h(\T_{\phi})=0$ (i.e. $f$ is maximal), then $\phi$ is of the
form $\phi=\e_ke_k^*$ with $\e_k=\pm1$. We observe that,
$|\phi(\sum_{k=1}^la_ke_k)|=|a_k|=a_k$. From \cite{AH} Proposition
4.8, we can choose
$\gamma\in\Gamma_{q_{k+1}-1}\backslash\Gamma_{q_k+k}$ such that
$|x_k(\gamma)|\geq\frac{1}{C}\|x_k\|=\frac{1}{C}$. Then,
$|\phi(\sum_{k=1}^la_ke_k)|=a_k\leq
C|a_k||x_k(\gamma)|=C|e_{\gamma}^*(a_kx_k)|\leq
C|e_{\gamma}^*(\sum_{k=1}^la_kx_k)|$.

We assume that for every $\phi\in W(\A_n,\bb)$ with
$h(\T_{\phi})=h>0$ and $\maxsupp\phi=l_0$, there exists
$\gamma\in\Gamma$, such that: \be \item
$\gamma\in\Gamma_{q_{l_0+1}+h}\backslash\Gamma_{q_{l_0+1}}$ \item
$h(\T_{\phi})=h(\F_{\gamma})\leq l_0$ \item
$|\phi(\sum_{k=1}^la_ke_k)|\leq C|\sum_{k=1}^la_kx_k(\gamma)|$ for
every $l\geq l_0$ \ee Observe that assumption (1) yields
$x_{l_0}<\rank\gamma<x_{l_0+1}$, while assumption (2) gives us that
$\minsupp x_{l_0+1}-\maxsupp x_{l_0}>h(\T_{\phi})$. Indeed,
\[x_{l_0}<q_{l_0+1}<\rank\gamma\leq q_{l_0+1}+h\leq
q_{l_0+1}+l_0<q_{l_0+1}+(l_0+1)<x_{l_0+1}\] \[\text{and } \minsupp
x_{l_0+1}-\maxsupp x_{l_0}>l_0+1>l_0\geq h(\F_{\gamma}).\]

Let $\phi\in W(\A_n,\bb)$ with $h(\T_{\phi})=h+1$,
$l_0=\maxsupp\phi$ and let $(\phi_t)_{t\in\T}$ the tree analysis
of $\phi$. Then, $\phi$ is of the form $\phi=\sum_{s\in
S_{\emptyset}}b_s\phi_s$, $\# S_{\emptyset}\leq n$. We observe
that for every $s\in S_{\emptyset}$, $h(\T_ {\phi_s})=h$. We set
$p_1=1$, for every $s\in S_{\emptyset}\diagdown\{1\}$
$p_s=\min\{q_k+k:k\in\supp\phi_s\}$ and for every $s\in S_{\emptyset}$, $r_s=q_{l_s+1}+h$ where $l_s=\maxsupp\phi_s$.\\
We next apply the inductive hypothesis to obtain
$\xi_s\in\Gamma_{r_s}\backslash\Gamma_{q_{l_s}+1}$ with
$h(\T_{\phi_s})=h(\F_{\xi_s})$ such that
\begin{eqnarray*}
|\phi_s(\sum_{k=1}^la_ke_k)| &=& |\phi_s(\sum_{k\in\supp
\phi_s}a_ke_k)|\leq C\e_s\sum_{k\in\supp
\phi_s}a_kx_k(\xi_s)\\ &=& C\e_se_{\xi_s}^*(\sum_{k\in\supp
\phi_s}a_kx_k) = C\e_se_{\xi_s}^*\circ
P_{(p_s,r_s]}(\sum_{k=1}^la_kx_k),
\end{eqnarray*}
with $\e_s$ such that
$\e_se_{\xi_s}^*(\sum_{k\in\supp\phi_s}a_kx_k)=|\sum_{k\in\supp
\phi_s}a_kx_k(\xi_s)|$.

Let $\gamma\in\Gamma$ have analysis
$\{p_s,r_s,\e_se_{\xi_s}^*\}_{s\in
S_{\emptyset}}\bigcup\{\eta_s\}_{s\in
S_{\emptyset}\diagdown\{1\}}$ where $\eta_s\in\Delta_{r_s+1}$.
Observe that $\rank\xi_s\in{(q_{l_s+1},r_s]}\subset{(p_s,r_s]}$.
It is clear that for every $s\in S_{\emptyset}\diagdown\{1\}$,
$d_{\eta_s}^*(\sum_{k=1}^la_kx_k)=0$. Indeed, \[\supp
x_{l_s}<q_{l_s+1}<q_{l_s+1}+{(h+1)}=r_s+1\leq
q_{l_s+1}+(l_s+1)<\supp x_{l_s+1}.\] Therefore,
\begin{eqnarray*}
|\phi(\sum_{k=1}^la_ke_k)| &\leq&\sum_{s\in
S_{\emptyset}}|b_s\phi_s(\sum_{k\in\supp\phi_s}a_ke_k)| \leq
C\sum_{s\in S_{\emptyset}}b_s\e_se_{\xi_s}^*\circ
P_{(p_s,r_s]}(\sum_{k=1}^la_kx_k)\\ &\leq&
C|\sum_{k=1}^la_kx_k(\gamma)|
\end{eqnarray*}
It is clear that $h(\T_{\phi})=h(\F_{\gamma})\leq l_0$ and
$x_{l_0}<\rank\gamma<x_{l_0+1}$.
 \ep

\bc\label{9} For every block sequence in $\mathfrak{X}_r$ there
exists a further block sequence satisfying inequality (\ref{eq1}).
\ec

%-------------------------Upper Inequality---------------------------%

\section{The upper estimate}
Let $(y_l)_{l\in\N}$ be a normalized skipped block sequence in
$\mathfrak{X}_r$. From Corollary \ref{9}, we can find a further block sequence of $(y_l)_l$,
still denoted by $(y_l)_l$, satisfying inequality (\ref{eq1}). \\
Therefore, we have that
\[\|\sum_{l=1}^my_l\|_{\infty}\geq\frac{1}{C}\|\sum_{l=1}^me_l\|_{\T(\A_n,\bb)}\]
For every $j\in\N$, set $M_j=\{1,2,\ldots,n\}^j$. It is easily
checked, after identifying $M_j$ with $\{1, \dots, n^j\}$ for every $j$, that the functional $f_j=\sum_{s\in
M_j}(\prod_{i=1}^jb_{s_i})e_s^*$ belongs to $W(\A_n,\bb)$ where
$s_i$ is the $i$-th coordinate of $s$, for each $i=1,2,\ldots,n$
and $\sum_{s\in M_j}\prod_{i=1}^jb_{s_i}=(\sum_{i=1}^nb_i)^j$.
Using the fact that $\#M_j=n^j$, we obtain that\\
$\|\sum_{l=1}^{n^j}e_l\|_{\T(\A_n,\bb)}=\|\sum_{s\in
M_j}e_s\|_{\T(\A_n,\bb)}\geq
f_j(\sum_{l=1}^{n^j}e_l)=(\sum_{i=1}^nb_i)^j$.

Also, for every $m\in\N$ large enough we may find $j\in\N$ such
that $n^{j+1}>m\geq n^j$. From the above and the unconditionality
of the basis of the space $\T(\A_n,\bb)$, it follows that
\[\|\sum_{l=1}^my_l\|_{\infty}\geq\frac{1}{C}\|\sum_{l=1}^me_l\|_{\T(\A_n,\bb)}\geq\frac{1}{C}\|\sum_{l=1}^{n^j}e_l\|_{\T(\A_n,\bb)}=(\sum_{i=1}^nb_i)^j\] We conclude that
$\|\sum_{l=1}^my_l\|_{\infty}\stackrel{m\to\infty}\longrightarrow\infty$ as $\sum_{i=1}^n b_i>1$.\\

We next choose a further block sequence $(x_k)_{k\in\N}$ of
$(y_l)_{l\in\N}$ with some additional properties. Let $\e>0$ and
choose a descending sequence $(\e_k)_k$ of positive reals such
that $(\sum_{k=1}^{\infty}\e_k)<\e$. We can also find an
increasing sequence $(n_k)_k$ of positive integers and a sequence
$(F_k)_k$ of succesive subsets of $\N$ such that the following are
satisfied:
\begin{enumerate}
\item For every $k\in\N$, $\frac{1}{n_k}<\e_k$.
\item For every $k\in\N$, $\|\sum_{l\in F_k}y_l\|> n_k$. This is possible, due to the above notation.
\end{enumerate}
We have thus constructed a normalized skipped block sequence $(x_k)_{k\in\N}$
of the form $x_k=\sum_{l\in F_k}\lambda_l y_l$, where $\lambda_l=\frac{1}{\|\sum_{l\in F_k}y_l\|} $. Notice that $|\lambda_l|<\e_k$ for every $l\in F_k$.

Let $\gamma\in\Gamma$ with tree analysis $\F_{\gamma}=
\{\xi_t,(p_t,q_t]\}_{t\in\T}$.\\
For every $k\in\N$, we set $t_k=\max\{t:\ran
x_k\subset(p_t,q_t]\}$. Notice that if for a given $x_k$, $t_k$ is
non-maximal, then there exist at least two immediate successors of
$t_k$, say $s_1$, $s_2$ such that the corresponding intervals
$(p_{s_1},q_{s_1}]$, $(p_{s_2},q_{s_2}]$ intersect $\ran x_k$. For
later use we shall denote by $(p_{s_0},q_{s_0}]$ the first
interval in the natural order of disjoint segments  of the natural
numbers that intersects
$x_k$. Notice that $s_0$ is not necessarily the first element of $S_t$.\\
For the pair $\gamma$, $(x_k)_{k\in\N}$ and for every $t\in\T$ we
define the following sets: $D_t=\bigcup_{s\succeq t}\{k: s=t_k\}$,
$K_t=D_t\backslash\cup_{s\in S_t}D_s\ =\{k:t=t_k\}$ and
$E_t=\{s\in S_t: D_s\neq\emptyset\}$.

We now set $x_k=x'_k+x''_k+x'''_k$ where,
\[x'_k=x_k\mid_{(p_{s_0},q_{s_0}]},\ \ x''_k=x_k\mid_{\bigcup_{s\in
S_{t_k},s\neq s_0}(p_s,q_s]}\ \text{and } x'''_k=x_k-x'_k-x''_k.\]

\begin{remark}\label{7}
\be \item The sets $D_t$,$K_t$,$E_t$ are determined by the chosen
pair $\gamma, (x_k)_k$. For a different pair, these sets may
differ as well. For example, let $k\in K_t$, for the pair $\gamma,
(x_k)_k$. Then $t=t_k$ for $x_k$. By the construction of $x'_k$,
there exists $s_k\in S_t$ such that
$x'_k=x_k\mid_{(p_{s_k},q_{s_k}]}$. Thus, taking the pair $\gamma,
(x'_k)_k$ the same $k$ belongs to $K_{s_k}$. \item For every
$k\in\N$,
$|g_{t_k}(x_k)|\leq 2Cn\e_k$.\\
Indeed, from the definition of $(x_k)_{k\in\N}$ we have that

\begin{eqnarray*}
|g_{t_k}(x_k)| &\leq& \sum_{s\in
S_{t_k}^{p_{t_k}}}|d_{\eta_s}^*(x_k)|\leq\sum_{s\in
S_{t_k}^{p_{t_k}}}|e_{\eta_s}^*\circ
P_{\{q_s+1\}}(\sum_{l\in F_k}\lambda_l y_l)|\leq\\
 &\leq& \sum_{s\in
S_{t_k}^{p_{t_k}}}\|e_{\eta_s}^*\|\|P_{\{q_s+1\}}\||\lambda_l^s|\|y_l^s\|
\leq\sum_{s\in S_{t_k}^{p_{t_k}}}2C\e_k\leq \\
&\leq& 2C\e_k(\sharp S_{t_k}) \leq 2Cn\e_k.
\end{eqnarray*}

\item It is obvious that $g_{t_k}(x_k)=g_{t_k}(x'''_k)$,
$f_{t_k}(x'''_k)=0$ and for every $t\prec t_k$, $g_t(x'''_k)=0$.
\ee
\end{remark}

\bl\label{6} For the pairs $\gamma, (x'_k)_{k\in\N}$ and $\gamma,
(x''_k)_{k\in\N}$ it holds that $\# (K_t\cup E_t)\leq n$. \el

\bp
Let $t\in\T$ and let $k\in K_t$.\\
We set $s_k=\max \{s\in S_t:(p_s,q_s]\cap\ran
x'_k\neq\emptyset\}$. From the definition of $t_k$, notice that
$\# S_t\geq2$. It holds that $s_k\not\in E_t$.\\
Indeed, from the definition of $t_k$, $s_k$ we have that
$(p_{t_k},q_{t_k}]\cap\ran x'_k=\ran x'_k$ and
$(p_{s_k},q_{s_k}]\cap\ran x'_k=(p_{s_k},q_{s_k}]$. Since $s_k\in
S_{t_k}$, $(p_{s_k},q_{s_k}]\subseteq(p_{t_k},q_{t_k}]$. It
follows that $(p_{s_k},q_{s_k}]\subseteq\ran x'_k$.\\ Therefore,
we can define a one-to-one map $G: K_t\rightarrow S_t\backslash
E_t$, hence $\# K_t
+\# E_t\leq \# S_t\leq n$.\\
The proof for the pair $\gamma, (x''_k)_{k\in\N}$ is similar. \ep

\bpr\label{8} Let $(x_k)_{k\in\N}$ be as above.Then for every
$\gamma\in\Gamma$ there exist $\phi_1,\phi_2\in W(\A_n,\bb)$ such
that for every sequence $(a_k)_{k\in\N}$ of positive scalars, for
every $l\in\N$ it holds that,
\begin{equation}\label{eq2}
|\sum_{k=1}^l a_kx_k(\gamma)|\leq
\frac{1}{b_n}(\phi_1+\phi_2)(\sum_{k=1}^l
a_ke_k)+2Cn\e(\sum_{k=1}^la_k^r)^{\frac{1}{r}}
\end{equation}
\epr

\bp Let $\gamma\in\Delta_{q+1} with a(\gamma)=a\leq n$. Let
$\F_{\gamma}= \{\xi_t,(p_t,q_t]\}_{t\in\T}$,
where $\xi_{\emptyset}=\gamma$, be the tree analysis of $\gamma$. We
may assume that $\bigcup_{k=1}^l\ran x_k\subset
(p_{\emptyset},q_{\emptyset}]$.

\begin{claim}
For the pairs $\gamma, (x'_k)_{k\in\N}$ and $\gamma,
(x''_k)_{k\in\N}$ there exist $\phi_1,\phi_2\in W(\A_n,\bb)$ such
that for every sequence of positive scalars $(a_k)_{k\in\N}$ and for
every $l\in\N$, it holds that
\begin{equation}\label{eq3}
|f_{\emptyset}(\sum_{k=1}^la_kx'_k)|\leq \frac{2C}{b_n}
\phi_1(\sum_{k=1}^la_ke_k)
\end{equation}
\begin{equation}\label{eq4}
|f_{\emptyset}(\sum_{k=1}^la_kx''_k)|\leq \frac{2C}{b_n}
\phi_2(\sum_{k=1}^la_ke_k)
\end{equation}
\end{claim}

\bp[Proof of the Claim]
We only prove inequality \ref{eq3}. The proof of inequality \ref{eq4} requires
the same arguments. We recall that $f_t=\sum_{s\in S_t}b_s\e_s(f_s+g_s)\circ P_{(p_s,q_s]}$ for every $t\in\T$ non maximal. From the definition of $(x'_k)_{k\in\N}$, we have that $g_s\circ P_{(p_s,q_s]}(x'_k)=0$ for every $s\in S_t$.Therefore,\\
 $f_t(\sum_{k\in D_t}a_kx'_k)=(\sum_{s\in S_t}b_s\e_sf_s\circ P_{(p_s,q_s]})(\sum_{k\in D_t}a_kx'_k)$. We will use backwards induction on the levels of the tree $\T$, i.e we shall
 show that for every $t\in\T$ there exists $\phi_1^t\in W(\A_n,\bb)$ with $\supp\phi_1^t\subseteq D_t$ such that \[|f_t(\sum_{k\in D_t}a_kx'_k)|\leq \frac{2C}{b_n}\phi_1^t(\sum_{k\in
     D_t}a_ke_k)\].\\
 Let $0<h\leq\max\{|t|:t\in\T\}$\\
 We assume that the proposition has been proved for all $t$ with $|t|=h$.\\
 Let $t\in\T$ with $|t|=h-1$.Then we have the following cases:
 \begin{enumerate}
 \item If $f_t$ is a maximal node, $f_t(\sum_{k\in D_t}a_kx'_k)=0$, so there is nothing to prove. Indeed, $K=D_t$, therefore for every $k\in D_t$, from Corollary \ref{5} $f_t(x'_k)=0$ since $t=t_k$.
 \item If $f_t$ is a non-maximal node, then\\ $f_t(\sum_{k\in D_t}a_kx'_k)=(\sum_{s\in S_t}b_s\e_sf_s\circ P_{(p_s,q_s]})(\sum_{k\in D_t}a_kx'_k)=$\\
     $=\sum_{s\in S_t}b_s\e_sf_s(\sum_{k\in D_s}a_kx'_k)+\sum_{k\in K}(\sum_{s\in S_t}b_s\e_sf_s)(a_kx'_k)$.\\
     From the fact that, for every $k\in K_t$, $g_t(x'_k)=0$ we get that
     \[|f_t(x'_k)|=|x'_k(\xi_t)|\leq\|x'_k\|\leq
     2C=2Ce_k^*(e_k).\]
     Moreover, for $s\in E_t$ it holds that $|s|=h-1$. For every $k\in D_s$, from the inductive hypothesis we obtain
     \[|\sum_{s\in
     S_t}b_sf_s(x'_k)|=|b_sf_s(x'_k)|\leq b_s\frac{2C}{b_n}\phi_1^s(e_k).\]
     with $\phi_1^s\in W(\A_n,\bb)$ and $\supp\phi_1^s\subseteq
     D_s$.\\
     We set $\phi_1^t=(\sum_{s\in E_t}b_s\phi_1^s+\sum_{k\in K_t}b_ke_k^*)$.\\
     From Lemma \ref{6}, it is easily checked that $\phi_1^t\in W(\A_n,\bb)$ and it holds that, $|f_t(\sum_{k\in D_t}a_kx'_k)|\leq \frac{2C}{b_n}\phi_1^t(\sum_{k\in
     D_t}a_ke_k)$.
 \end{enumerate}
\ep
Recall that
$e_{\gamma}^*(\sum_{k=1}^la_kx_k)=g_{\emptyset}(\sum_{k=1}^la_kx_k)+f_{\emptyset}(\sum_{k=1}^la_kx_k)$.\\
The fact that
$g_{\emptyset}(\sum_{k=1}^la_kx'_k)=g_{\emptyset}(\sum_{k=1}^la_kx''_k)=g_{\emptyset}(\sum_{k\in{\{m:t_m\neq\emptyset\}}}a_kx'''_k)=f_{\emptyset}(\sum_{k\in{\{m:t_m=\emptyset\}}}a_kx'''_k)=0$
implies the following:
\begin{eqnarray*}
|e_{\gamma}^*(\sum_{k=1}^la_kx_k)| &\leq&
|g_{\emptyset}(\sum_{k\in{\{m:t_m=\emptyset\}}}a_kx'''_k)|+|f_{\emptyset}(\sum_{k=1}^la_kx'_k)|\\
&+&
|f_{\emptyset}(\sum_{k=1}^la_kx''_k)|+|f_{\emptyset}(\sum_{k\in{\{m:t_m\neq\emptyset\}}}a_kx'''_k)|
\end{eqnarray*}
>From Remark \ref{7} we get that,
\[|g_{\emptyset}(\sum_{k\in{\{m:t_m=\emptyset\}}}a_kx'''_k)|\leq
\sum_{k\in{\{m:t_m=\emptyset\}}}a_k|g_{\emptyset}(x'''_k)| \leq
2Cn\sum_{k\in{\{m:t_m=\emptyset\}}}a_k\e_k.\] From Lemma \ref{4}
and Remark \ref{7} we have that,
\begin{eqnarray*}
|f_{\emptyset}(\sum_{k\in{\{m:t_m\neq\emptyset\}}}a_kx'''_k)|
&\leq&
\sum_{k\in{\{m:t_m\neq\emptyset\}}}a_k(\prod_{t<t_k}b_t)|g_{t_k}(x'''_k)|\leq \\ &\leq&
2C\frac{1}{2}n\sum_{k\in{\{m:t_m\neq\emptyset\}}}a_k\e_k
\leq 2Cn\sum_{k\in{\{m:t_m\neq\emptyset\}}}a_k\e_k.
\end{eqnarray*}
 Finally, we conclude that
\begin{eqnarray*}
|\sum_{k=1}^la_kx_k(\gamma)| &\leq&
2Cn\sum_{k\in{\{m:t_m=\emptyset\}}}a_k\e_k+\frac{2C}{b_n}\phi_1(\sum_{k=1}^la_ke_k)\\
&+&
\frac{2C}{b_n}\phi_2(\sum_{k=1}^la_ke_k)+2Cn\sum_{k\in{\{m:t_m\neq\emptyset\}}}a_k\e_k\\
&\leq& \frac{2C}{b_n}(\phi_1+\phi_2)(\sum_{k=1}^l
a_ke_k)+2Cn\sum_{k=1}^la_k\e_k\\ &\leq&
\frac{2C}{b_n}(\phi_1+\phi_2)(\sum_{k=1}^l
a_ke_k)+2Cn\max\{a_k:k\in\N\}(\sum_{k=1}^l\e_k)\\
&\leq& \frac{2C}{b_n}(\phi_1+\phi_2)(\sum_{k=1}^l
a_ke_k)+2Cn\e(\sum_{k=1}^la_k^r)^{\frac{1}{r}}.
\end{eqnarray*}
 where in the last inequality we used the fact that the $\ell_r$
 norm dominates the $c_0$ norm.
\ep

\br\label{11} From \cite{AT} Theorem I.4, we know that $\|\sum
a_ke_k\|_{\T(\A_n,\bb)}\geq M(\sum a_k^r)^{\frac{1}{r}}$. This
result and the previous Proposition, yield that
\[|\sum_{k=1}^l
a_kx_k(\gamma)|\leq\frac{2C}{b_n}(\phi_1+\phi_2)(\sum_{k=1}^l
a_ke_k)+\frac{2Cn\e}{M}\|\sum_{k=1}^la_ke_k\|_{\T(\A_n,\bb)}.\] For
$\e=\frac{M}{n b_n}$,
\[|\sum_{k=1}^l
a_kx_k(\gamma)|\leq\frac{6C}{b_n}\|\sum_{k=1}^l
a_ke_k\|_{\T(\A_n,\bb)}.\]
Therefore,
\begin{equation}\label{eq5}
\|\sum_{k=1}^la_kx_k\|_{\infty}\leq\frac{6C}{b_n}\|\sum_{k=1}^l
a_ke_k\|_{\T(\A_n,\bb)}.
\end{equation}
\er

\bc\label{10} For every block sequence in $\mathfrak{X}_r$ there
exists a further block sequence satisfying inequality (\ref{eq5}).
\ec

%-------------------------Main Result---------------------------%

\section{The main result}

\bpr\label{12} Let $(x_k)_{k\in\N}$ be a skipped block sequence in
$\mathfrak{X}_r$ satisfying $\minsupp x_{k+1}>\maxsupp x_k+k$ and
the conditions of Proposition \ref{8}. Then $(x_k)_{k\in\N}$ is
equivalent to the basis of the Tsirelson space $\T(\A_n,\bb)$ for
$n$ and $\bb$ determined as before. \epr

\bp It is an immediate consequence  of Propositions \ref{3},
\ref{8} and Remark \ref{11}. \ep

\bpr \label{13} The space $\T(\A_n,\bb)$ is isomorphic to $\ell_p$
for some $p\in(1,\infty)$. \epr

\bp In a similar manner as in \cite{AT} Theorem I.4, one can see
that for every normalized block sequence $(x_k)_k$ of the basis
$(e_j)_j$ and for every scalar sequence $(a_k)$ it holds that, $\|\sum
a_kx_k\|\leq\frac{2}{b_n}\|\sum a_ke_k\|$. Zippin's Theorem
\cite{Z} yields that $\T(\A_n,\bb)$ is isomorphic to some $\ell_p$
for some $p\in(1,\infty)$.\ep

\br An alternative proof could also be derived using the Results in
Sections 4 and 5. Indeed, let $(y_l)_{l\in\N}$ be a skipped block
sequence in $\mathfrak{X}_r$. Then, there exists a further block
sequence $(x_k)_{k\in\N}$ satisfying simultaneously the
assumptions of Corollaries \ref{9} and \ref{10}. Therefore,
$(x_k)_{k\in\N}$ satisfies the assumptions of Proposition
\ref{12}.

Let's observe that every further block sequence
$(z_k)_k$ of $(x_k)_k$ is also skipped block and satisfies
Proposition \ref{12}, thus it is equivalent to the basis of the
space $\T(\A_n,\bb)$. Hence, every block sequence $(z_n)_n$ of
$(x_k)_k$ is equivalent to $(x_k)_k$. Zippin's theorem \cite{Z}
yields that the space $\overline{<(x_k)_k>}$ is isomorphic to some
$\ell_p$. Therefore, $\T(\A_n,\bb)\cong\ell_p$ for some
$p\in(1,\infty)$.\er

In order to determine the exact value of $p$, we need the
following Proposition.

\bpr \label{14} The space $\T(\A_n,\bb)$ is isomorphic to $\ell_r$
with $\frac{1}{r}+\frac{1}{r'}=1$ and $\sum_{i=1}^nb_i^{r'}=1$.
\epr

\bp First, let observe that for every $x\in c_{00}$,
$\|x\|\leq\|x\|_r$. We shall use induction on the cardinality of
$\supp x$. If $|\supp x|=1$, it is trivial. Assume that it holds
for every $y\in c_{00}$ with $|\supp y|\leq n$ and let $x\in
c_{00}$ with $|\supp x|=n+1$. Then either $\|x\|=\|x\|_{\infty}$
or $\|x\|=\sum_{i=1}^nb_i\|E_ix\|$ for some appropriate subsets
$E_1<E_2<\ldots<E_n$. In the first case, there is nothing to prove
as for every $p\in[r,\infty)$ $\|x\|_{\infty}\leq\|x\|_p$.
Therefore we only need to deal with the second case.

It suffices
to observe that for every $i=1,2,\ldots,n$, the cardinality of
$\supp E_ix$ is less than $\supp x$ and thus, using the inductive
hypothesis along with $H\ddot{o}lder's$ inequality, we get that
\[\|x\|\leq\sum_{i=1}^nb_i\|E_ix\|_r\leq(\sum_{i=1}^nb_i^{r'})^{\frac{1}{r'}}(\sum_{i=1}^n\|E_ix\|_r^r)^{\frac{1}{r}}=\|x\|_r.\]
By combining
the preceding argument with Proposition \ref{13}, we conclude that
$\T(\A_n,\bb)$ is isomorphic to $\ell_p$ for some
$p\in[r,\infty)$.\\ For every $l\in\N$ set
$M_l=\{1,2,\ldots,n\}^l$. We have already mentioned that for every
$l\in\N$ the functional $f_l=\sum_{s\in
M_l}(\prod_{i=1}^lb_{s_i})e_s^*$ belongs to $W(\A_n,\bb)$ where
$s_i$ is the $i$-th coordinate of $s$, for each $i=1,2,\ldots,n$
and $\sum_{s\in M_l}\prod_{i=1}^lb_{s_i}=(\sum_{i=1}^nb_i)^l$. We
set $a_s=\prod_{i=1}^lb_{s_i}$ and $x_l=\sum_{s\in
M_l}a_s^{\frac{r'}{r}}e_s$. It is easily seen that for every
$l\in\N$, $\|x_l\|=1$. Indeed, \[\|x_l\|\leq\|x_l\|_r=(\sum_{s\in
M_l}a_s^{r'})^{\frac{1}{r}}=(\sum_{i=1}^nb_i^{r'})^{\frac{l}{r}}=1=f_l(x_l)\leq\|x_l\|.\]

We claim that for $p'>r$ and every $\e>0$ there exists $l\in\N$
such that $\|x_l\|_{p'}<\e$.
If the claim holds we are done as $p$ coincides with $r$.\\
\textit{Proof of the Claim}: Notice that for $p'>r$,
$\sum_{i=1}^nb_i^{\frac{r'}{r}p'}=\sum_{i=1}^nb_i^{r'(1+\delta)}$
for some $0<\delta<1$. But for every $i=1,2,\ldots,n$ $b_i<1$,
and therefore
\[\sum_{i=1}^nb_i^{r'(1+\delta)}<\sum_{i=1}^nb_i^{r'}=1.\]\\
Thus,
there exists $l\in\N$ such that
$(\sum_{i=1}^nb_i^{r'(1+\delta)})^l<\e^{p'}$. Then for this $l$,
\[\|x_l\|_{p'}=(\sum_{s\in
M_l}a_s^{\frac{r'}{r}p'})^{\frac{1}{p'}}=(\sum_{s\in
M_l}a_s^{r'(1+\delta)})^{\frac{1}{p'}}=(\sum_{i=1}^nb_i^{r'(1+\delta)})^{\frac{l}{p'}}<\e.\]
\ep

\bt For every $r\in(1,\infty)$ the space $\mathfrak{X}_r$ is
$\ell_r$ saturated. \et

\bp As it was mentioned in the above Remark, for every skipped
block sequence in $\mathfrak{X}_r$ we can find a further block
sequence $(x_k)_k$ such that the space $\overline{<(x_k)_k>}$ is
isomorphic to $\ell_r$. \ep

\br From the previous Theorem, we deduce that the space
$\mathfrak{X}_r$ is a separable $\mathcal{L}^{\infty}$ space which
does not contain $\ell_1$. Therefore, the results of
D.Lewis-C.Stegall \cite{LS} and A. Pelczy\'nski \cite{P} yields
that $\mathfrak{X}_r^*$ is isomorphic to $\ell_1$. Alternatively, one
can use the corresponding argument of D. Alspach \cite{Al} and show
directly that $(M_q)$ is a shrinking FDD for $\mathfrak{X}_r$. It then follows
that $(e_\gamma^*)_{\gamma \in \Gamma}$ is a basis for $\mathfrak{X}_r^*$, equivalent
to the usual $\ell_1$-basis.
\er

%-------------------------Bibliography---------------------------%

\end{document}